\documentclass[12pt,leqno]{article}
\usepackage{amsthm}
\usepackage{amsmath}
\usepackage{amssymb}
\usepackage{amsfonts}
\usepackage[center]{titlesec}

\begin{document}
\newtheoremstyle{mytheoremstyle} 
    {\topsep}                    
    {\topsep}                    
    {\itshape}                   
    {}                           
    {\scshape}                   
    {.}                          
    {.5em}                       
    {}  
\theoremstyle{mytheoremstyle}    
\newcounter{thmcounter}
\newcounter{Remarkcounter}
\newcounter{examcon}
\numberwithin{thmcounter}{section}
\numberwithin{examcon}{section}
\newtheorem{Remark}[Remarkcounter]{Remark}
\newtheorem{Prop}[thmcounter]{PROPOSITION}
\newtheorem{Corol}[thmcounter]{COROLLARY}
\newtheorem{theorem}[thmcounter]{THEOREM}
\newtheorem{Lemma}[thmcounter]{LEMMA}
\newtheorem{example}[examcon]{EXAMPLE}
\newcommand{\tr}{\mathrm{tr}}

\title{\textbf{ON ROTATIONAL SURFACES WITH ZERO MEAN CURVATURE IN THE PSEUDO--EUCLIDEAN SPACE $\mathbb{E}_2^4$}}
\author{B. BEKTA\c{S}$^1$, E. \"{O}. CANFES$^2$ and U. DURSUN$^3$}

\date{}
\maketitle

{\scriptsize{$^{1,2}$ Faculty of Science and Letters, Department of Mathematics, Istanbul Technical University, 34469,
Istanbul, Turkey, \centerline{e-mail: bektasbu@itu.edu.tr, canfes@itu.edu.tr}}}

{\scriptsize{$^{3}$ Department of Mathematics, I\c{s}{\i}k University, \c{S}ile Campus, 34980, \c{S}ile, Istanbul, Turkey  
\\ \centerline{e-mail: ugur.dursun@isikun.edu.tr}}}

\renewcommand{\thefootnote}{}

\textbf{Abstract.}
In this work, we study a class of rotational surfaces in the pseudo--Euclidean space $\mathbb{E}_2^4$ whose profile curves lie in two--dimensional planes. We solve the differential equation that characterizes the rotational surfaces with zero mean curvature to determine the profile curves of such rotational surfaces. Then, we give some explicit parametrization of maximal rotational surfaces and the timelike surfaces with zero mean curvature in $\mathbb{E}_2^4$.

\section{Introduction}
Minimal surfaces are important in geometry as well as in physics. Minimal surfaces in real space forms and indefinite space forms have been studied extensively by many mathematicians. In \cite{Moore}, Moore studied general rotational surfaces in the four dimensional Euclidean space $\mathbb{E}^4$. He proved that if there exists a minimal general rotational surface in $\mathbb{E}^4$ with equal rates of rotation, then its profile curves lie in 2--planes, and he determined such minimal surfaces. In \cite{Arslan etc}, some special solutions of the differential equation that characterizes minimal general rotational surfaces in $\mathbb{E}^4$ with profile curves lie in 2--planes were given. Then, the third author and  Turgay obtained the general solution of the differential equation that characterizes minimal rotational surfaces with different rates of rotation in $\mathbb{E}^4$, \cite{DursunTurgay}.

On the other hand, Ganchev and Milousheva used some special invariant to classify  minimal general rotational surfaces in the Euclidean space $\mathbb{E}^4$ and the Minkowski space $\mathbb{E}_1^4$, \cite{GanchMilous3,GanchMilous2}. They also obtained all timelike surfaces in $\mathbb{E}_1^4$ with zero mean curvature in the class of rotational surfaces of Moore type, \cite{GanchMilous1}. In \cite{HuiLiGuiLi}, the particular solutions of differential equation related to minimal surfaces in the pseudo--Euclidean space $\mathbb{E}_2^4$ were obtained. Recently, Chen studied minimal Lorentz surfaces in arbitrary indefinite space forms. In \cite{Chen}, he obtained several classification results, in particular, he completely classified all minimal Lorentz surfaces with arbitrary dimension $m$ and arbitrary index $s$.  
 
In this paper, we study rotational surfaces in the pseudo--Euclidean space $\mathbb{E}^4_2$ with profile curves lying in 2--planes. We solve the differential equation that characterizes the rotational surfaces in $\mathbb{E}_2^4$ with zero mean curvature. Thus, we give examples of maximal surfaces and Lorentz surfaces with zero mean curvature in $\mathbb{E}_2^4$.

\section{Prelimineries}
Let $\mathbb E^m_s$ be the $m$--dimensional pseudo--Euclidean space with the canonical
metric given by
$$
 \tilde{g}=\sum_{i=1}^{m-s}(dx_i)^2-\sum_{i=m-s+1}^m(dx_i)^2, 
$$
where $(x_1, x_2, \ldots, x_m)$ is a standard rectangular coordinate system in $\mathbb E^m_s$.

For a point $c\in\mathbb{E}_s^m$ and $r>0$, we put
\begin{eqnarray}
\nonumber
\mathbb{S}^{m-1}_{s}(c, r)&=&\{x\in\mathbb{E}_s^m|\langle x-c, x-c\rangle =r^2\},\\
\mathbb{H}^{m-1}_{s-1}(c, -r)&=&\{x\in\mathbb{E}_s^m|\langle x-c, x-c\rangle =-r^2\}\notag,
\end{eqnarray} 
where $\langle , \rangle$ denotes the indefinite inner product on $\mathbb{E}_s^m$. $\mathbb{S}^{m-1}_{s}(c, r)$ and $\mathbb{H}^{m-1}_{s-1}(c, -r)$ are called, respectively, a pseudo--sphere and a pseudo--hyperbolic space. The hyperbolic space $\mathbb{H}^{m-1}(c, -r)$ is defined by
\begin{equation}
\nonumber
\mathbb{H}^{m-1}(c, -r)=\{x\in\mathbb{E}_1^m|\langle x-c, x-c\rangle =-r^2, x_m>0\}.
\end{equation} 
If $c=0$, then $\mathbb{S}^{m-1}_{s}(0, r)$ and $\mathbb{H}^{m-1}_{s-1}(0, -r)$ are denoted by $\mathbb{S}^{m-1}_{s}(r)$ and $\mathbb{H}^{m-1}_{s-1}(-r)$.

A vector $v \in \mathbb{E}^m_s$ is called spacelike (resp., timelike)
if $\langle v, v\rangle>0$ or $v=0$ (resp., $\langle v,v\rangle<0$). A vector $v$ is called lightlike if $\langle v,v\rangle=0$ and $v\neq 0$.  A submanifold $M$ of $\mathbb{E}^{m}_s$ is said to be spacelike if every non--zero tangent vector on $M$ is spacelike and it is said to be timelike if at least one of non--zero tangent vector on $M$ is timelike. 

Let $M$ be an oriented $n$--dimensional submanifold in the $(n+2)$--dimensional
pseudo--Euclidean space $\mathbb E^{n+2}_2$. We choose an oriented local orthonormal  frame $\{e_1, \dots, $ $ e_{n+2} \}$ on $M$ with $\varepsilon_A= \langle e_A, e_A \rangle = \pm 1$ such that  $e_1,\dots,e_{n}$ are tangent to $M$ and $e_{n+1}, \,e_{n+2}$ are normal to $M$. We use the following convention on the range of indices: $1\leq i, j, k, \ldots\leq n$, $n+1\leq r, s, t, \ldots\leq n+2$.
 
Let $\widetilde{\nabla}$ be the Levi--Civita connection of $\mathbb E^{n+2}_2$
and  $\nabla$  the induced connection on $M$. Denote by $\{\omega^1,\dots,\omega^{n+2} \}\, $ the dual frame  and by  $\{\omega_{AB}\},  A,B=1,\dots, n+2,$  the
connection forms associated to $\{e_1,\dots,e_{n+2} \}$. The formulas of Gauss and Weingarten are given by, respectively,
\begin{align}
\nonumber
\widetilde{\nabla}_{e_k}e_i= \sum_{j=1}^{n}\varepsilon_j \omega_{ij}(e_k) e_j + \sum_{r=n+1}^{n+2} \varepsilon_r h^r_{ik} e_r, \;\;\mbox{and}\;\;
\widetilde{\nabla}_{e_k}e_r=  - A_r(e_k)+ \sum_{s=n+1}^{n+2} \varepsilon_s \omega_{rs}(e_k) e_s, 
\end{align}
where $h^r_{ij}$ is the coefficients of the second fundamental form $h$, and
$A_r$ the Weingarten map in the direction $e_r$.

The mean curvature vector $H$ is defined by $H= \frac{1}{n}\sum_{r,i} \varepsilon_i\varepsilon_r h^r_{ii} e_r$. A spacelike submanifold $M$ in $\mathbb{E}^m_s$ is called maximal if $H$ vanishes identically.  

The Codazzi equation of $M$ in $\mathbb{E}^{n+2}_2$ is given by
\begin{align} \label{codazzi}
\begin{split}
    &h^r_{ij,k}= h^r_{jk,i},\\
  &h^r_{jk,i} = e_i(h^r_{jk}) + \sum_{s=n+1}^{n+2} \varepsilon_s h^s_{jk} \omega_{sr} (e_i) -
\sum_{\ell=1}^{n} \varepsilon_{\ell}\left ( \omega_{j\ell}(e_i) h^r_{\ell k} +
\omega_{k \ell}(e_i) h^r_{\ell j} \right ).
\end{split}
\end{align}

 
\subsection{Rotational surfaces in $\mathbb{E}_2^4$}
Let $M_1(b)$ and $M_2(b)$ be rotational surfaces in the pseudo--Euclidean space $\mathbb{E}_2^4$ whose profile curves lie in 2--planes. 
We can choose a profile curve $\alpha$ of $M_1(b)$ in the $yw$--plane as $\alpha(u)=(0, y(u), 0, w(u))$, defined on an open interval $I\subset\mathbb{R}$, and thus the parametrization of $M_1(b)$ is given by
\begin{equation}\label{1DY}
M_1(b): r_1(u,v)=(w(u)\sinh v, y(u)\cosh (bv), y(u)\sinh (bv), w(u)\cosh v), 
\end{equation}
for some constant $b>0$, where $u\in I$ and $v\in\mathbb{R}$.

We consider the following orthonormal moving frame field $\{e_1, e_2, e_3, e_4\}$ on $M_1(b)$ such that $e_1, e_2$ are tangent to $M_1(b)$, and $e_3, e_4$ are normal to $M_1(b)$:
\begin{eqnarray}
\label{TegNormVektors1DY} e_1&=&\frac{1}{q}\frac{\partial}{\partial v}, \quad
e_2 =\frac{1}{A}\frac{\partial}{\partial u},\\
\label{TegNormVektors1DY3} e_3&=&\frac{1}{A}(y'(u)\sinh v, w'(u)\cosh (bv), w'(u)\sinh (bv), y'(u)\cosh v),\\
\label{TegNormVektors1DY4} e_4&=&-\frac{\varepsilon\varepsilon^*}{q}(by(u)\cosh v, w(u)\sinh (bv), w(u)\cosh (bv), by(u)\sinh v),
\end{eqnarray}
where $A=\sqrt{\varepsilon({y^\prime}^2(u)-{w^\prime}^2(u))}\neq 0$, $q= \sqrt{\varepsilon^*(w^2(u)- b^2 y^2(u))}\neq 0$, and 
$\varepsilon=\mbox{sgn}{({y^\prime}^2(u)-{w^\prime}^2(u))}$, $\varepsilon^*=\mbox{sgn}(w^2(u)- b^2 y^2(u))$.
Then, $\varepsilon_1 = -\varepsilon_4=\varepsilon^* , \; \varepsilon_2 = -\varepsilon_3 =\varepsilon$.

By a direct calculation, we have the components of the second fundamental form and the connection forms as follows
\begin{align}
\label{1DYuzeyTemelForm3}
h^3_{11}&=\frac{1}{Aq^2}(b^2y(u)w'(u)-w(u)y'(u)),\quad h^3_{22}=\frac{1}{A^3}(w'(u)y''(u)-y'(u)w''(u)),   \\
\label{1DYuzeyTemelForm4}
h^4_{12}&=\frac{\varepsilon\varepsilon^* b}{Aq^2}(w(u)y'(u)-y(u)w'(u)), \quad h^3_{12}=h^4_{11}=h^4_{22}=0, \\
\label{1DYuzeyTegetKonForm}
\omega_{12}(e_1)&= \frac{1}{Aq^2}(b^2y(u)y'(u)-w(u)w'(u)), \quad \omega_{12}(e_2)=0, \\
\label{1DYuzeyNormlKonForm}
\omega_{34}(e_1)&= \frac{\varepsilon\varepsilon^* b}{Aq^2}(w(u)w'(u)-y(u)y'(u)),\quad \omega_{34}(e_2)=0.
\end{align}

We can choose a profile curve $\beta$ of $M_2(b)$ in the $xz$--plane as $\beta(u)=(x(u), 0, z(u), 0)$ defined on an open interval $I\subset\mathbb{R}$, and thus the parametrization of $M_2(b)$ is given by
\begin{equation}\label{2DY}
M_2(b): r_2(u,v)=(x(u)\cos v, x(u)\sin v, z(u)\cos (bv), z(u)\sin (bv)) 
\end{equation}
for some constant $b>0$, where $u\in I$ and $v\in (0, 2\pi)$.

We consider the following orthonormal moving frame field $\{e_1, e_2, e_3, e_4\}$ on $M_2(b)$ such that $e_1, e_2$ are tangent to $M_2(b)$, and $e_3, e_4$ are normal to $M_2(b)$:
\begin{eqnarray}
\label{TegNormVektors2DY} e_1&=&\frac{1}{\bar{q}}\frac{\partial}{\partial v}, \quad
e_2 =\frac{1}{\bar{A}}\frac{\partial}{\partial u},\\
\label{TegNormVektors2DY3} e_3&=&\frac{1}{\bar{A}}(z'(u)\cos v, z'(u)\sin v, x'(u)\cos (bv), x'(u)\sin (bv)),\\
\label{TegNormVektors2DY4} e_4&=&-\frac{\varepsilon\varepsilon^*}{\bar{q}}(bz(u)\sin v, -bz(u)\cos v, x(u)\sin (bv), -x(u)\cos (bv)),
\end{eqnarray}
where $\bar{A}=\sqrt{\varepsilon({x^\prime}^2(u)-{z^\prime}^2(u))}\neq 0$, $\bar{q}= \sqrt{\varepsilon^*(x^2(u)- b^2 z^2(u))}\neq 0$, and $\varepsilon=\mbox{sgn}{({x^\prime}^2(u)-{z^\prime}^2(u))}$, $\varepsilon^*=\mbox{sgn}(x^2(u)- b^2 z^2(u))$.
Then, $\varepsilon_1 = -\varepsilon_4=\varepsilon^* , \; \varepsilon_2 = -\varepsilon_3 =\varepsilon$.

By a direct computation, we have the components of the second fundamental form and the connection forms as follows
\begin{align}
\label{2DYuzeyTemelForm3}
h^3_{11}&=\frac{1}{\bar{A}{\bar{q}}^2}(b^2z(u)x'(u)-x(u)z'(u)),\quad h^3_{22}=\frac{1}{{\bar{A}}^3}(z'(u)x''(u)-x'(u)z''(u)),   \\
\label{2DYuzeyTemelForm4}
h^4_{12}&=\frac{\varepsilon\varepsilon^* b}{\bar{A}{\bar{q}}^2}(z(u)x'(u)-x(u)z'(u)), \quad h^3_{12}=h^4_{11}=h^4_{22}=0, \\
\label{2DYuzeyTegetKonForm}
\omega_{12}(e_1)&= \frac{1}{\bar{A}{\bar{q}}^2}(b^2z(u)z'(u)-x(u)x'(u)), \quad \omega_{12}(e_2)=0, \\
\label{2DYuzeyNormlKonForm}
\omega_{34}(e_1)&= \frac{\varepsilon\varepsilon^* b}{\bar{A}{\bar{q}}^2}(z(u)z'(u)-x(u)x'(u)),\quad \omega_{34}(e_2)=0.
\end{align}
Therefore, we have the mean curvature vector for the rotational surfaces $M_1(b)$ and $M_2(b)$ as follows
\begin{eqnarray}
\label{12DYOrtEgrVek} H&=&-\frac 12(\varepsilon\varepsilon^* h^3_{11} + h^3_{22})e_3 .
\end{eqnarray}
On the other hand, by using the Codazzi equation \eqref{codazzi} we obtain
\begin{eqnarray}
\label{12DYUzCodazzi1} e_2(h^3_{11})&=&\varepsilon^* h_{12}^4\omega_{34}(e_1)+\omega_{12}(e_1)(\varepsilon^* h_{11}^3-\varepsilon h_{22}^3),\\
\label{12DYUzCodazzi2}e_2(h^4_{12})&=&-\varepsilon h_{22}^3\omega_{34}(e_1) + 2 \varepsilon^*h_{12}^4\omega_{12}(e_1).
\end{eqnarray}
The rotational surfaces $M_1(b)$ and $M_2(b)$ defined by \eqref{1DY} and \eqref{2DY} for $b=1$, $x(u)=y(u)=f(u)\sinh u$ and $z(u)=w(u)=f(u)\cosh u$ are also known as Vranceanu rotational surface, where $f(u)$ is a smooth function, \cite{HuiLiGuiLi}.  

\section{Rotational Surfaces with Zero Mean Curvature}
In this section, we determine all rotational surfaces $M_1(b)$ and $M_2(b)$ defined, respectively, by \eqref{1DY} and \eqref{2DY} with zero mean curvature.

By considering \eqref{1DYuzeyTemelForm3} and \eqref{12DYOrtEgrVek}, a rotational surface $M_1(b)$ has zero mean curvature if and only if the coordinate functions $y(u)$ and $w(u)$ of the profile curve $\alpha$ satisfy the differential equation
\begin{equation}
\label{difdenk1}
w'(u)y''(u)-y'(u)w''(u)+
({y^\prime}^2(u)-{w^\prime}^2(u))\frac{b^2y(u)w'(u)-w(u)y'(u)}{w^2(u)-b^2y^2(u)}=0.
\end{equation}
Note that $y(u)=c_0w(u), c_0^2\neq 1$, is a solution of differential equation \eqref{difdenk1} for $b=1$, and it can be shown easily that $M_1(1)$ is an open part of a timelike plane in $\mathbb{E}_2^4$. Thus, we rule out this case. 

\begin{Prop}
\label{propM1zeromeancurvatureb1}
A non--planar rotational surface $M_1(b)$ in $\mathbb{E}_2^4$ defined by \eqref{1DY}  for $b=1$ has zero mean curvature if and only if its profile curve is given by
\begin{equation}
\label{sol1}
(y(u)+w(u))^2+\lambda_0(w(u)-y(u))^2=\mu_0
\end{equation}
for some constants $\lambda_0\neq 0$ and $\mu_0$. 
\end{Prop}

PROOF.
Assume that $M_1(1)$ has zero mean curvature. So, for $b=1$ the differential equation \eqref{difdenk1} can be written as
\begin{equation}
\nonumber
\frac{\left(\frac{y'(u)}{w'(u)}\right)^\prime}{1-\left(\frac{y'(u)}{w'(u)}\right)^2}+\frac{\left(\frac{w(u)}{y(u)}\right)^\prime}{1-\left(\frac{w(u)}{y(u)}\right)^2}=0
\end{equation}
from which the first integration gives
$\tanh^{-1}\left(\frac{y'(u)}{w'(u)}\right)+\tanh^{-1}\left(\frac{w(u)}{y(u)}\right)=c$ for some constant $c$. Using the logarithmic expression for the inverse hyperbolic tangent function we get 
\begin{equation}
\nonumber
(y(u)+w(u))(y'(u)+w'(u))\pm e^{2c}(y(u)-w(u))(y'(u)-w'(u))=0.
\end{equation}
The solution of this differential equation yields \eqref{sol1} for some constant $\lambda_0=\pm e^{2c}\neq 0$ and $\mu_0$. 

The converse of the proof of the theorem comes from direct computation. 
\qed

The solution \eqref{sol1} is a quadratic curve. For some suitable values of $\lambda_0$ and $\mu_0$, we have ellipses or hyperbolas. For instance, if we take $\lambda_0=1$ and $\mu_0=2$, then from \eqref{sol1} we have $w^2(u)+y^2(u)=1$, that is, the profile curve $\alpha$ is a part of the unit circle, $w^2+y^2=1$. When we choose $y(u)=\sin u$ and $w(u)=\cos u$, we have $\varepsilon=\varepsilon^*=\mbox{sgn}({\cos 2u})$. Hence, the surface $M_1(1)$ is maximal for $|u|<\frac{\pi}{4}$, and the parametrization of $M_1(b)$ becomes 
\begin{equation}
\nonumber
M_1(1): r_1(u, v)=(\cos u\sinh v, \sin u\cosh v, \sin u\sinh v, \cos u\cosh v)
\end{equation}
for $u\in(-\frac{\pi}{4},\frac{\pi}{4})$ and $v\in\mathbb{R}$. Similarly, we can take $y(u)=\cos u$ and $w(u)=\sin u$. In this case, we have  $\varepsilon=\varepsilon^*=-\mbox{sgn}({\cos 2u})$, and by choosing $\frac{\pi}{4}<u<\frac{3\pi}{4}$, the surface $M_1(1)$ is maximal with positive definite metric, and by choosing $|u|<\frac{\pi}{4}$, $M_1(1)$ is a maximal surface with negative definite metric.

If we take $\lambda_0=-1$ and $\mu_0=4$, then from \eqref{sol1} we get $y(u)w(u)=1$, that is, the profile curve $\alpha$ is the part of the hyperbola, $yw=1$. By taking $y(u)=u$ and $\displaystyle{w(u)=\frac{1}{u}, u>0}$, the parametrization of $M_1(b)$ is given by
\begin{equation}
\nonumber
M_1(1): r_1(u, v)=\left(\frac{\sinh v}{u}, u\cosh v, u\sinh v, \frac{\cosh v}{u}\right)
\end{equation}
which is regular for $0<u<1$ or $u>1$, and it is timelike with zero mean curvature.
 
In \cite{HuiLiGuiLi}, it was shown that the Vranceanu rotational surface is maximal if $f(u)=a(\cosh(2u+c))^{-1/2}$, where $a$ and $c$ are constants. Also, for this function $f(u)$, the component functions $y(u)$ and $w(u)$ satisfy the equation \eqref{sol1}.

We assume that the profile curves $\alpha$ and $\beta$ of $M_1(b)$ and $M_2(b)$, respectively, are arc length parametrized, that is,
${y^\prime}^2(u)-{\omega^\prime}^2(u)=\varepsilon$ and ${x^\prime}^2(u)-{z^\prime}^2(u)=\varepsilon$.

Now, we give the following lemma to obtain the general solution of the differential equation \eqref{difdenk1} for $b\neq 1$: 
\begin{Lemma}\label{LemmaM1ZeroMeanCurvature}
Let $M_1(b)$ be a non--planar rotational surface in the pseudo--Euclidean space $\mathbb{E}^4_2$ given by \eqref{1DY} with $b\neq 1$. Then, $M_1(b)$ has zero mean curvature if and only if the component functions $y(u)$ and $w(u)$ of the unit speed profile curve $\alpha$ of $M_1(b)$ satisfy 
the differential equation 
\begin{equation}
\label{M1denk0}
(b^2-1)(b^2y^2(u){w'}^2(u)-w^2(u){y'}^2(u))=a_0
\end{equation}
for some constant $a_0$ and an open subinterval $J\subset I$ on which $y'(u)w'(u)\ne 0$. 
\end{Lemma}

PROOF.
Let $M_1(b)$ be a non--planar rotational surface in $\mathbb{E}_2^4$ defined by \eqref{1DY} with $b\neq 1$ such that the profile curve $\alpha$ of $M_1(b)$ is unit speed. Assume that the mean curvature of $M_1(b)$ is zero, i.e., $\varepsilon\varepsilon^* h_{11}^3=-h_{22}^3$. If $y=y_0=constant$ or $w=w_0=constant$ on an open subinterval of $I$, then $M_1(b)$ is a planar rotational surface or has no zero mean curvature. So there is an open subinterval $J\subset I$ on which $y'(u)w'(u)\neq 0$. Using $\varepsilon\varepsilon^* h_{11}^3=-h_{22}^3$ in the Codazzi equations \eqref{12DYUzCodazzi1} and \eqref{12DYUzCodazzi2}, we obtain that
\begin{equation}
\label{M1denk1}
e_2(h^3_{11})=\varepsilon^*h^4_{12}\omega_{34}(e_1)+2\varepsilon^*h^3_{11}\omega_{12}(e_1),
\end{equation}
\begin{equation}
\label{M1denk2}
e_2(h^4_{12})=\varepsilon^*h^3_{11}\omega_{34}(e_1)+2\varepsilon^*h^4_{12}\omega_{12}(e_1).
\end{equation}
From \eqref{M1denk1} and \eqref{M1denk2}, we get
\begin{equation}
\label{M1denk3}
h^3_{11}e_2(h^3_{11})-h^4_{12}e_2(h^4_{12})
=2\varepsilon^*\left((h^3_{11})^2-(h^4_{12})^2\right)\omega_{12}(e_1).
\end{equation}
It is clear that $(h_{11}^3)^2-(h_{12}^4)^2=0$ is a solution of \eqref{M1denk3}. In this case, by considering the first equations in \eqref{1DYuzeyTemelForm3} and \eqref{1DYuzeyTemelForm4}, and $b\neq 1$, we have \eqref{M1denk0} for  $a_0=0$.

On the other hand, if $(h^3_{11})^2-(h^4_{12})^2\ne 0$ on some open subinterval $J\subset I$, by using \eqref{1DYuzeyTegetKonForm} and \eqref{M1denk3} we then obtain that
\begin{equation}
\frac{e_2((h_{11}^3)^2-(h_{12}^4)^2)}{4((h_{11}^3)^2-(h_{12}^4)^2)}+\frac{w(u)w'(u)-b^2y(u)y'(u)}{w^2(u)-b^2y^2(u)}=0.
\end{equation}
By integrating this equation, we get
\begin{equation}
\label{M1denk5}
((h^3_{11})^2-(h^4_{12})^2)(w^2(u)-b^2y^2(u))^2=a_0
\end{equation}
for some constant $a_0\neq 0$.
Hence, using the first equations in \eqref{1DYuzeyTemelForm3} and \eqref{1DYuzeyTemelForm4}, equation \eqref{M1denk5} yields \eqref{M1denk0}.

Conversely, assume that the coordinate functions $y(u)$ and $w(u)$ of the unit speed profile curve $\alpha$ satisfy 
the differential equation \eqref{M1denk0} for some constant $a_0$, and $y'(u)w'(u)\neq 0$ on an open subinterval $J\subset I$. From ${y^\prime}^2(u)-{w^\prime}^2(u)=\varepsilon$ and \eqref{M1denk0}, we get
\begin{equation}
\label{M1denk6}
{y^\prime}^2(u)=\frac{\tilde{a}_0+\varepsilon b^2y^2(u)}{b^2y^2(u)-w^2(u)}\;\;\mbox{and}\;\; {w^\prime}^2(u)=\frac{\tilde{a}_0+\varepsilon w^2(u)}{b^2y^2(u)-w^2(u)},  
\end{equation}
where $\displaystyle{{\tilde a}_0}=\frac{a_0}{b^2-1}$.
Differentiating equations in \eqref{M1denk6} with respect to $u$ and using again \eqref{M1denk6}, we obtain 
\begin{equation}
\nonumber
y'(u)y''(u)=\frac{\varepsilon b^2y(u)y'(u)}{b^2y^2(u)-w^2(u)}-{y^\prime}^2(u)\frac{b^2y(u)y'(u)-w(u)w'(u)}{b^2y^2(u)-w^2(u)},
\end{equation}
\begin{equation}
\nonumber
w'(u)w''(u)=\frac{\varepsilon w(u)w'(u)}{b^2y^2(u)-w^2(u)}-{w^\prime}^2(u)\frac{b^2y(u)y'(u)-w(u)w'(u)}{b^2y^2(u)-w^2(u)}\cdot
\end{equation}
If we multiply these equations by $-{w^\prime}^2(u)$ and ${y^\prime}^2(u)$, respectively, and add them, we get 
\begin{equation}
y'(u)w'(u)\left(w'(u)y''(u)-y'(u)w''(u)+
\varepsilon\frac{b^2y(u)w'(u)-w(u)y'(u)}{w^2(u)-b^2y^2(u)}\right)=0
\end{equation}
which implies \eqref{difdenk1} as $y'(u)w'(u)\neq 0$ on $J\subset I$. Hence, the surface $M_1(b)$ has zero mean curvature.  
\qed

\begin{theorem}\label{TheoremM1ZeroMeanCurvature}
Let $M_1(b)$ be a non--planar rotational surface in the pseudo--Euclidean space $\mathbb{E}^4_2$ defined by \eqref{1DY} with $b\neq 1$. Then, $M_1(b)$ has zero mean curvature if and only if the component functions $y(u)$ and $w(u)$ of the unit speed profile curve $\alpha$ of $M_1(b)$ satisfy one of the following regular curves: 
\begin{itemize}
\item[i.] For $\varepsilon\varepsilon^*=1$,
\begin{equation}
\label{M1A0(2)}
\sin^{-1}\left(\frac{w(u)}{\mu_0}\right)=\pm\frac{1}{b}\sin^{-1}\left(\frac{by(u)}{\mu_0}\right)+c_0,\;\;\;\; \mu_0= \sqrt{\frac{\varepsilon^*a_0}{1-b^2}},
\end{equation}
where $a_0\neq 0$ and $c_0$ are constants such that $\displaystyle{\frac{\varepsilon^*a_0}{1-b^2}>0}$. In this case, the surface $M_1(b)$ is spacelike with positive or negative definite metric.
\item[ii.] For $\varepsilon\varepsilon^*=-1$,
\begin{equation}
\label{M1A0(1)}
\left(w(u)+\sqrt{w^2(u)-\mu_0^2}\;\right)^{\pm{b}}=d_0\left(by(u)+\sqrt{b^2y^2(u)-\mu_0^2}\;\right),\;\;\;\; \mu_0=\sqrt{\frac{\varepsilon^*a_0}{b^2-1}},
\end{equation}
where $a_0$ and $d_0\neq 0$ are constants. If $a_0=0$, then $y(u)=b_0(w(u))^{\pm b}$, where $b_0$ is non--zero constant. In this case, the surface $M_1(b)$ is timelike.
\end{itemize}
\end{theorem}

PROOF.
Let $M_1(b)$ be a non--planar rotational surface in the pseudo--Euclidean space $\mathbb{E}_2^4$ given by \eqref{1DY} with $b\neq 1$ and zero mean curvature. Then, Lemma \ref{LemmaM1ZeroMeanCurvature} implies that the coordinate functions $y(u)$ and $w(u)$ of the unit speed profile curve $\alpha$ of $M_1(b)$ on an open subinterval of $I$ on which $y'(u)w'(u)\neq 0$ satisfy \eqref{M1denk0} for some constant $a_0$. 
Using ${y^\prime}^2(u)-{w^\prime}^2(u)=\varepsilon$ and \eqref{M1denk0}, we have \eqref{M1denk6} from which we get 
\begin{equation}
\label{M1denk7}
\sqrt{-\varepsilon^*(\tilde{a}_0+\varepsilon b^2y^2(u))}\;w'(u)=\pm\sqrt{-\varepsilon^*(\tilde{a}_0+\varepsilon w^2(u))}\;y'(u),
\end{equation}
where $\tilde{a}_0=\frac{a_0}{b^2-1}$. It has a solution according to $\varepsilon^*\varepsilon=1$ or $\varepsilon^*\varepsilon=-1$. \\
\textit{Case 1}: $\varepsilon\varepsilon^*=1$ and $a_0\neq 0$. Then, equation \eqref{M1denk7} becomes 
\begin{equation}
\label{M1denk7(1)}
\frac{dw}{\sqrt{-w^2(u)-\varepsilon^*{\tilde{a}}_0}}=\pm\frac{dy}{\sqrt{-b^2y^2(u)-\varepsilon^*{\tilde{a}_0}}}.
\end{equation}
If $\varepsilon^*\tilde{a}_0=\mu_0^2>0$, then $-(w^2(u)+\mu_0^2)<0$ and $-(b^2y^2(u)+\mu_0^2)<0$, and thus there is no solution of \eqref{M1denk7(1)}. Let $\varepsilon^*\tilde{a}_0=-\mu_0^2<0$, that is, $\mu_0^2=\frac{\varepsilon^*a_0}{1-b^2}>0$. The solution of equation \eqref{M1denk7(1)} is given by  
\begin{equation}
\nonumber
\sin^{-1}\left(\frac{w(u)}{\mu_0}\right)=\pm\frac{1}{b}\sin^{-1}\left(\frac{by(u)}{\mu_0}\right)+ c_0
\end{equation} 
for some constant $c_0$. In this case, the surface $M_1(b)$ is spacelike with positive or negative definite metric. \\
\textit{Case 2}: $\varepsilon\varepsilon^*=-1$. Then, equation \eqref{M1denk7} becomes 
\begin{equation}
\label{M1denk7(2)}
\frac{dw}{\sqrt{w^2(u)-\varepsilon^*\tilde{a}_0}}=\pm\frac{dy}{\sqrt{b^2y^2(u)-\varepsilon^*\tilde{a}_0}}.
\end{equation}
For $\varepsilon^*\tilde{a}_0=\mu_0^2>0$, the general solution of \eqref{M1denk7(2)} is given by 
\begin{equation}
\nonumber
\left(w(u)+\sqrt{w^2(u)-\mu_0^2\;}\right)^{\pm b}=d_0\left(by(u)+\sqrt{b^2y^2(u)-\mu_0^2\;}\right),
\end{equation} 
where $d_0$ is a non--zero constant. By a similar calculation, we can have \eqref{M1A0(1)} for $\varepsilon^*\tilde{a}_0=-\mu_0^2<0$. In this case, $M_1(b)$ is a timelike surface.

Conversely, let $M_1(b)$ be a rotational surface given by \eqref{1DY} whose profile curve $\alpha$ is given by one of the regular curves \eqref{M1A0(2)} and \eqref{M1A0(1)}. Suppose that the component functions of the profile curve $\alpha$ satisfy \eqref{M1A0(2)}. By differentiating \eqref{M1A0(2)} with respect to $u$, we obtain \eqref{M1denk7(1)} which yields $y'(u)w'(u)\neq 0$ on the interval $I$ since $\alpha$ is regular. From ${y^\prime}^2(u)-{w^\prime}^2(u)=\varepsilon$ and \eqref{M1denk7(1)}, we get
\begin{equation}
\nonumber
{y^\prime}^2(u)=\frac{\tilde{a}_0+\varepsilon b^2y^2(u)}{b^2y^2(u)-w^2(u)}\;\;\mbox{and}\;\; {w^\prime}^2(u)=\frac{\tilde{a}_0+\varepsilon w^2(u)}{b^2y^2(u)-w^2(u)}  
\end{equation}
which satisfy \eqref{M1denk0}. By a similar argument, it can be shown that the regular curve given by \eqref{M1A0(1)} satisfies \eqref{M1denk0}. Thus, by Lemma \ref{LemmaM1ZeroMeanCurvature}, the rotational surface $M_1(b)$ has zero mean curvature in the pseudo--Euclidean space $\mathbb{E}_2^4$.
\qed

Now we will give some examples of rotational surfaces with zero mean curvature whose profile curve $\alpha$ given by \eqref{M1A0(2)} or \eqref{M1A0(1)}.
\begin{example}
\upshape
For $\varepsilon=\varepsilon^*=1$, if we choose $a_0=\frac{3}{4}$, $b=\frac{1}{2}$ and $c_0=0$, then from \eqref{M1A0(2)} we obtain $\displaystyle{\sin^{-1}(w(u))=2\sin^{-1}\left(\frac{y(u)}{2}\right)}$. When we take $y(u)=2\sin u$, then we have $w(u)=\sin(2u)$. Thus, the parametrization of the surface $M_1(b)$ becomes
\begin{equation}
\nonumber
M_1(1/2): r_1(u, v)=\left(\sin(2u)\sinh v, 2\sin u\cosh\left(\frac{v}{2}\right), 2\sin u\sinh\left(\frac{v}{2}\right), \sin(2u)\cosh v\right)
\end{equation}
which is a maximal surface with positive definite metric for $0<u<\frac{\pi}{4}$ and $v\in\mathbb{R}$.
\end{example}
\begin{example}
\upshape 
For $\varepsilon=-\varepsilon^*=1$, if we choose $a_0=-3$, $b=2$ and $d_0=1$, then from \eqref{M1A0(1)} we can take the component functions of the profile curve $\alpha$ as $y(u)=\frac{1}{2}\cosh(2u)$ and $w(u)=\cosh u$. Then, the parametrization of  $M_1(b)$ becomes
\begin{equation}
\nonumber
M_1(2): r_1(u, v)=\left(\cosh u\sinh v, \frac{1}{2}\cosh(2u)\cosh(2v), \frac{1}{2}\cosh(2u)\sinh(2v), \cosh u\cosh v\right)
\end{equation}
which is a timelike surface with zero mean curvature in $\mathbb{E}_2^4$ for $u>0$ and $v\in\mathbb{R}$.\\
\end{example}
\begin{example}
\upshape
If we choose $b_0=1$ and $b=2$, we have $y(u)=(w(u))^2$ from the equation \eqref{M1A0(1)} for constant $a_0=0$. Let $y(u)=u^2$ and $w(u)=u$, $u>0$ be the parametrization of $y=w^2$. Thus, the parametrization of the surface $M_1(b)$ is given by
\begin{equation}
\nonumber
M_1(2): r_1(u,v)=\left(u\sinh v, u^2\cosh(2v), u^2\sinh(2v), u\cosh v\right)
\end{equation}
which is a timelike surface with zero mean curvature for $0<u<\frac{1}{2}$ or $u>\frac{1}{2}$,  and $v\in\mathbb{R}$.
\end{example}
By a similar way, we study the rotational surface $M_2(b)$ given by \eqref{2DY} with zero mean curvature. Considering \eqref{2DYuzeyTemelForm3} and \eqref{12DYOrtEgrVek}, a rotational surface $M_2(b)$ has zero mean curvature if and only if the coordinate functions $x(u)$ and $z(u)$ of the profile curve $\beta$ satisfy the differential equation
\begin{equation}
\label{2difdenk1}
z'(u)x''(u)-x'(u)z''(u)+
({x^\prime}^2(u)-{z^\prime}^2(u))\frac{b^2z(u)x'(u)-x(u)z'(u)}{x^2(u)-b^2z^2(u)}=0.
\end{equation}
Note that $x(u)=cz(u), c^2\neq 1$, is a solution of differential equation \eqref{2difdenk1} for $b=1$. But in this case it can be shown easily that $M_2(1)$ is an open part of a spacelike plane in $\mathbb{E}_2^4$. Thus, we rule out this case.

We state the following proposition for the solution of \eqref{2difdenk1} without proof because its proof is similar to the proof of Proposition \ref{propM1zeromeancurvatureb1}.  
\begin{Prop}
\label{propM2zeromeancurvatureb1}
A non--planar rotational surface $M_2(b)$ in $\mathbb{E}_2^4$ defined by \eqref{2DY} for $b=1$ has zero mean curvature if and only if its profile curve is given by
\begin{equation}
\label{2sol1}
(x(u)+z(u))^2+\lambda_0(x(u)-z(u))^2=\mu_0
\end{equation}
for some constants $\lambda_0\neq 0$ and $\mu_0$. 
\end{Prop}
The solution \eqref{2sol1} is a quadratic curve. For some suitable values of $\lambda_0$ and $\mu_0$, we have ellipses or hyperbolas. For instance, if we take $\lambda_0=1$ and $\mu_0=2$, then we have $x^2(u)+z^2(u)=1$ from \eqref{2sol1}, that is, the profile curve $\beta$ is a part of the unit circle, $x^2+z^2=1$. When we choose $x(u)=\cos u$ and $z(u)=\sin u$, we have $\varepsilon^*=-\varepsilon=\mbox{sgn}(\cos 2u)$. Hence, the surface $M_2(1)$ is timelike in $\mathbb{E}_2^4$, and its parametrization is given by
\begin{equation} 
\nonumber
M_2(1): r_2(u,v)=(\cos u\cos v, \cos u\sin v, \sin u\cos v, \sin u\sin v)
\end{equation}
for $u\in(-\frac{\pi}{4},\frac{\pi}{4})$ and $v\in(0, 2\pi)$.
Similarly, if we can take $x(u)=\sin u$ and $z(u)=\cos u$, then the surface $M_2(1)$ is again timelike with zero mean curvature for $|u|<\frac{\pi}{4}$.

If we take $\lambda_0=-1$ and $\mu_0=4$, then from \eqref{2sol1} we get $x(u)z(u)=1$,
that is, the profile curve $\beta$ is the part of the hyperbola $xz=1$. By taking $x(u)=u$ and $\displaystyle{z(u)=\frac{1}{u}, u>0}$, the parametrization of $M_2(b)$ is given by 
\begin{equation}
\nonumber
M_2(1): r_2(u,v)=\left(u\cos v, u\sin v, \frac{\cos v}{u}, \frac{\sin v}{u}\right)
\end{equation}
which is maximal with positive or negative definite metric according to $u\!>1$ or $0<u<1$, respectively. 

In \cite{HuiLiGuiLi}, it was shown that the Vranceanu rotational surface has zero mean curvature 
if $f(u)=a(\cosh(2u+c))^{{-1}/{2}}$, where $a$ and $c$ are constants and it is timelike with zero mean curvature. Also, for this function $f(u)$, the component functions $x(u)$ and $z(u)$ satisfies the equation \eqref{2sol1}.

The formulas for $M_1(b)$ such as second fundamental form and differential equation of zero mean curvature are valid for $M_2(b)$ if we replace $y(u)$ and $w(u)$ with $z(u)$ and $x(u)$, respectively. For that replacement only the sign of ${x^\prime}^2(u)-{z^\prime}^2(u)$ changes, that is, ${y^\prime}^2(u)-{w^\prime}^2(u)=\varepsilon$ turns to be ${x^\prime}^2(u)-{z^\prime}^2(u)=-\varepsilon$. 

Thus, we give the following lemma and theorem without proof because their proofs are similar to the proof of Lemma \ref{LemmaM1ZeroMeanCurvature} and Theorem \ref{TheoremM1ZeroMeanCurvature}.

\begin{Lemma}\label{LemmaM2ZeroMeanCurvature}
Let $M_2(b)$ be a non--planar rotational surface in the pseudo--Euclidean space $\mathbb{E}^4_2$ given by \eqref{2DY} 
with $b\neq 1$. Then, $M_2(b)$ has zero mean curvature if and only if the coordinate functions $x(u)$ and $z(u)$ of the unit speed profile curve $\beta$ of $M_2(b)$ satisfy 
the differential equation 
\begin{equation}
\label{M2denk0}
(b^2-1)(b^2z^2(u){x^\prime}^2(u)-x^2(u){z^\prime}^2(u))=\bar{a}_0
\end{equation}
for some constant $\bar{a}_0$ and an open interval $J\subset I$ on which $x'(u)z'(u)\ne 0$. 
\end{Lemma}

\begin{theorem}\label{TheoremM2ZeroMeanCurvature}
Let $M_2(b)$ be a non--planar rotational surface in the pseudo-Euclidean space $\mathbb{E}^4_2$ given by \eqref{2DY} with 
$b\neq 1$. Then, $M_2(b)$ has zero mean curvature if and only if the coordinate functions $x(u)$ and $z(u)$ of the unit speed profile curve $\beta$ of $M_2(b)$ satisfy one of the following regular curves:
\begin{itemize}
\item[i.] For $\varepsilon\varepsilon^*=1$, 
\begin{equation}\label{M21}
\left(x(u)+\sqrt{x^2(u)-\bar{\mu}_0^2}\;\right)^{\pm{b}}=\bar{c}_0\left (bz(u)+\sqrt{b^2z^2(u)-\bar{\mu}_0^2}\;\right),\;\;\;\; \bar{\mu}_0=\sqrt{\frac{\varepsilon^*\bar{a}_0}{b^2-1}},
\end{equation}
where $\bar{a}_0$ and $\bar{c}_0\neq 0$ are constants. If $\bar{a}_0=0$, then $z(u)=\bar{b}_0(x(u))^{\pm{b}}$ where $\bar{b}_0$ is non--zero constant. In this case, the surface $M_2(b)$ is spacelike with positive or negative definite metric.
\item[ii.] For $\varepsilon\varepsilon^*=-1$,
\begin{equation}\label{M21-}
\sin^{-1}\left(\frac{x(u)}{\bar{\mu}_0}\right)=\pm\frac{1}{b}\sin^{-1}\left(\frac{bz(u)}{\bar{\mu}_0}\;\right)+\bar{d}_0, \;\;\;\; \bar{\mu}_0=\sqrt{\frac{\varepsilon^*\bar{a}_0}{1-b^2}},
\end{equation} 
where $\bar{a}_0\neq 0$ and $\bar{d}_0$ are constants such that $\displaystyle{\frac{\varepsilon^*\bar{a}_0}{1-b^2}>0}$. In this case, the surface $M_2(b)$ is timelike. 
\end{itemize}  
\end{theorem}
Now, we will give some parametrization for the rotational surface $M_2(b)$ whose the profile curve $\beta$ is given by \eqref{M21} or \eqref{M21-}.
\begin{example}
\upshape
For $\varepsilon=\varepsilon^*=1$, if we choose $\bar{a}_0=3$, $b=2$ and $\bar{c}_0=e$, then from \eqref{M21} we can take the component functions of the profile curve $\beta$ as $x(u)=\cosh u$ and $z(u)=\frac{1}{2}\cosh(2u-1)$. Then, the rotational surface $M_2(b)$ defined by
\begin{equation}
\nonumber
M_2(2): r_2(u, v)=\left(\cosh u\cos v, \cosh u\sin v, \frac{1}{2}\cosh(2u-1)\cos(2v), \frac{1}{2}\cosh(2u-1)\sin(2v)\right)\\
\end{equation}
is maximal in $\mathbb{E}_2^4$ for $0<u<1$ and $v\in(0, 2\pi)$.
\end{example}
\begin{example}
\upshape
If we choose $\bar{b}_0=1$ and $b=2$, we have $z(u)=(x(u))^2$ from the equation \eqref{M21} for constant $\bar{a}_0=0$. Let $x(u)=u$ and $z(u)=u^2$, $u>0$ be the parametrization of $z=x^2$. Thus, the parametrization of $M_2(b)$ is given by
\begin{equation}
\nonumber
M_2(2): r_2(u,v)=\left(u\cos v, u\sin v, u^2\cos(2v), u^2\sin(2v)\right)
\end{equation}
which is a maximal surface with positive or negative definite metric, respectively, for $0<u<\frac{1}{2}$ or $u>\frac{1}{2}$.
\end{example}
\begin{example}
\upshape
For $\varepsilon=-\varepsilon^*=1$, if we choose $\bar{a}_0=-\frac{3}{4}$, $b=\frac{1}{2}$ and $\bar{d}_0=-\frac{\pi}{4}$, then from \eqref{M21-} we obtain that $\displaystyle{\sin^{-1}(x(u))=2\sin^{-1}\left(\frac{z(u)}{2}\right)-\frac{\pi}{4}}$. When we take $z(u)=2\sin u$, then we have $x(u)=\sin(2u-\frac{\pi}{4})$. Thus, the parametrization of $M_2(b)$ becomes
\begin{equation}
\nonumber
M_2(1/2): r_2(u, v)=\left(\sin\left(2u-\frac{\pi}{4}\right)\cos v, \sin\left(2u-\frac{\pi}{4}\right)\sin v, 2\sin u\cos\left(\frac{v}{2}\right), 2\sin u\sin\left(\frac{v}{2}\right)\right)\\
\end{equation}
which is a timelike surface with zero mean curvature for $\frac{\pi}{8}<u<\frac{\pi}{4}$ and $v\in(0, 2\pi)$.
\end{example}

\end{document}